\baselineskip=15pt \magnification=\magstep1 \font\msbm=msbm10
\leftskip1.5cm \hsize13cm

\centerline{\bf ON SOME FINITE SUMS WITH FACTORIALS} \vskip1cm
\centerline{\bf Branko Dragovich} \vskip1cm

\noindent{\bf Abstract.} The summation formula
$$
\sum^{n-1}_{i=0}\varepsilon^i\, i!\, (i^k+u_k) =
v_k+\varepsilon^{n-1}\, n!\, A_{k-1}(n)
$$
$(\varepsilon=\pm 1;\ k=1,2,\cdots;\ u_k,\, v_k\in
\msbm\hbox{Z};\, \,  A_{k-1}$ is a polynomial ) is derived and its
various aspects are considered. In particular, divisibility with
respect to $n$ is investigated. Infinitely many equivalents to
Kurepa's hypothesis on the left factorial are found.

\vskip1cm
\centerline{{\bf 1. Introduction}}
\vskip.5cm

The subject of the present paper is the investigation of  finite sums
of the form
$$
\sum^{n-1}_{i=0}\, \varepsilon^i\, i!\, P_k(i)\ ,\eqno(1)
$$
where $\varepsilon = \pm 1$, and
$$
P_k(i) = C_k\, i^k+\cdots+C_1\, i+C_0         \eqno(2)
$$
is a polynomial with $k,i\in\msbm\hbox{N}_0 =
\msbm\hbox{N}\cup\{0\}$ and coefficients $C_0,\, C_1,\cdots,C_k
\in\msbm\hbox{Z}$.

We mainly consider the following three problems of (1): a)
summation formula, b) divisibility by $n!$ and c) connection with
the Kurepa hypothesis (KH) on the left factorial.\footnote{}{{\it
1991 Mathematics Subject Classification.} Primary 11A05} All these
problems depend on the form of the polynomial $P_k(i)$ and  have
something in common with it.

In Sec. 2 we find a few ways to determine $P_k(i)$ which give
simple and useful summation formulae. Sec. 3 contains divisibility
properties.  The results concerning KH on the left factorial are
given in Sec. 4. Infinitely many equivalents to KH are found.

\vfill\eject \centerline{{\bf 2. Summation Formulae}} \vskip.5cm

\noindent{\bf Lemma 1.} {\it Let $\varepsilon = \pm 1$ and
$$
A_{k-1}(n) = a_{k-1}\, n^{k-1}+\cdots+a_1\, n+a_0\ ,\quad
k\in\msbm\hbox{N}\ , \quad n\in\msbm\hbox{N}_0\ ,\eqno(3)
$$
is a polynomial. One can find coefficients $a_{k-1} = 1$ and
$a_{k-2},\cdots,a_0\in\msbm\hbox{Z}$ such that identity
$$
(n+1)\, A_{k-1}(n+1)-\varepsilon\, A_{k-1}(n) =
n^k+A_{k-1}(1)-\varepsilon\, A_{k-1}(0)  \eqno(4)
$$
holds for all $n\in\msbm\hbox{N}_0$.}

 \vskip.5cm {\it Proof.}
Formula (4) has the form
$$
(n+1)\, A_{k-1}(n+1)-\varepsilon A_{k-1}(n) = n^k+u_k\ .\eqno(5)
$$
Replacing $A_{k-1}(n)$ by (3) and demanding (5) to be an identity, the
following system of linear equations must be satisfied:
$$
\eqalign{&\pmatrix{k\cr0\cr}a_{k-1} = 1\cr
&\bigg[\pmatrix{k\cr1\cr}-\varepsilon\bigg]a_{k-1}+a_{k-2} = 0\cr
&\pmatrix{k\cr2\cr}a_{k-1}+\bigg[\pmatrix{k-1\cr1\cr}-\varepsilon\bigg]a_{k-2}+
a_{k-3} = 0\cr
&\cdots\cdots\cdots\cdots\cdots\cdots\cdots\cdots\cdots\cdots\cdots\cdots\cdots\cr
&k\, a_{k-1}+(k-1)\, a_{k-2}+\cdots+(2-\varepsilon)\, a_1+a_0 =
0\cr &a_{k-1}+a_{k-2}+\cdots+a_1+(1-\varepsilon)\, a_0 = u_k\
.\cr}\eqno(6)
$$
Starting from the first equation, which gives $a_{k-1} = 1$, one
can in a successive way obtain solution for all $a_i = a_i(k,
\varepsilon), \ i = 0,\cdots,k-2$. The last equation in (6) serves
to determine $u_k$. Thus we get
$$
u_k = \sum^{k-1}_{i=0}a_i-\varepsilon \, a_0 =
A_{k-1}(1)-\varepsilon \, A_{k-1}(0)\ .\eqno(7)
$$

Note that (4) is an identity if and only if the  coefficients of the
polynomial $A_{k-1}(n)$ satisfy the system of linear equations (6),
where $u_k$ is given by (7).

The first five polynomials which satisfy (4) are:
$$
\eqalign{&A_0(n) = 1\ ,\cr &A_1(n) = n+\varepsilon-2\ ,\cr &A_2(n)
= n^2+(\varepsilon-3)\, n+4-5\, \varepsilon\ ,\cr &A_3(n) =
n^3+(\varepsilon-4)\, n^2+7\, (1-\varepsilon)\, n+18\,
\varepsilon-13\ ,\cr &A_4(n) = n^4+(\varepsilon-5)\, n^3+(11-9\,
\varepsilon)\, n^2+2\, (16\, \varepsilon-11)\, n+58-63\,
\varepsilon\ .\cr}\eqno(8)
$$
\vskip.5cm

\noindent{\bf Theorem 1.} {\it The summation formula
$$
\sum^{n-1}_{i=0} \varepsilon^i \, i! \,[i^k+A_{k-1}(1)-\varepsilon
\, A_{k-1}(0)] = -\varepsilon \, A_{k-1}(0)+\varepsilon^{n-1}\,
n!\, A_{k-1}(n)\eqno(9)
$$
is valid if and only if the polynomials $A_{k-1}(n), \
k\in\msbm\hbox{N}$, satisfy the identity (4).}

 \vskip.5cm

{\it Proof.} Summation of (4), previously multiplied by
$\varepsilon^i \, i!$, gives
$$
\eqalign{ & \sum^{n-1}_{i=0} \varepsilon^i \, i!\,
[i^k+A_{k-1}(1)-\varepsilon \, A_{k-1}(0)] \cr  = &
\sum^{n-1}_{i=0} \varepsilon^i \, i!\, [(i+1)\,
A_{k-1}(i+1)-\varepsilon A_{k-1}(i)] \, .\cr}  \eqno(10)
$$
Since on the r. h. s. all but the first and the last term cancel
we get (9). Now one can easily show that starting from (9) one
obtains (4).

Denoting $u_k = A_{k-1}(1)-\varepsilon \,  A_{k-1}(0),\, \ v_k =
-\varepsilon \, A_{k-1}(0)$ we can rewrite (9) in the form
$$
\sum^{n-1}_{i=0} \varepsilon^i \, i! \, (i^k+u_k) =
v_k+\varepsilon^{n-1}\, n!\, A_{k-1}(n)\ ,\quad k\ge1\ .\eqno(11)
$$

Formula (9), as well  as  (11), is determined by polynomial
$A_{k-1}(n)$ in (3), whose coefficients are solution of (6).
However, for large $k$, (6) becomes inconvenient. Therefore, it is
of interest to have another approach which is more effective to
get (11). \vskip.5cm

\noindent {\bf Theorem 2.} {\it If $\delta_{0k}$ is the Kronecker
symbol and
$$
S_k^\varepsilon(n) = \sum^{n-1}_{i=0} \varepsilon^i \, i!\, i^k\
,\quad\varepsilon = \pm 1\ ,\quad  k\in\msbm\hbox{N}_0\ ,\eqno(12)
$$
then
$$
S^\varepsilon_k(n) =
\delta_{0k}+\varepsilon\sum^{k+1}_{l=0}\pmatrix{k+1\cr
l\cr}S^\varepsilon_l(n)-\varepsilon^n \, n!\, n^k\ ,\quad
k\in\msbm\hbox{N}_0\ ,\eqno(13)
$$
is a recurrent relation.}

 \vskip.5cm {\it Proof.}
$$
\eqalign{S^\varepsilon_k(n) &= \delta_{0k}+\sum^{n-2}_{i=0}
\varepsilon^{i+1} \, (i+1)! \, (i+1)^k \cr & =
\delta_{0k}+\varepsilon \sum^{n-1}_{i=0}\varepsilon^i \, i! \,
(i+1)^{k+1}-\varepsilon^n \, n!\, n^k \cr & =
\delta_{0k}+\varepsilon \sum^{k+1}_{l=0}\pmatrix{k+1\cr l\cr} \,
S^\varepsilon_l(n)-\varepsilon^n \, n!\, n^k\ .\cr}
$$

Relation (13) gives a simpler way to find (11) in the explicit
form for a particular index $k \geq 0$.

From (13) one can obtain recurrent relations for $u_k$ and $v_k$.
In particular, when $\varepsilon = 1$, we have
$$
u_{k+1} = -k\, u_k-\sum^{k-1}_{l=1}\pmatrix{k+1\cr l\cr} \, u_l+1\
,\quad u_1 = 0\ ,\ \ k\ge 1\ ,\eqno(13.a)
$$
$$
v_{k+1} = -k\, v_k-\sum^{k-1}_{l=1}\pmatrix{k+1\cr l\cr} \,
v_l-\delta_{0k}\ ,\quad k\ge 0\ .\eqno(13.b)
$$

Some first values of $u_k$ and $v_k $  $(\varepsilon = 1)$ are:

\def\prored{&\vbox to 4pt{}&&&&&&&&&&&&&&&&&&&&&\cr}%
$$
\vbox{\offinterlineskip\halign{\tabskip=4pt#\vrule&\hfil$#$\hfil&\vrule#&\hfil$#$\hfil&\vrule#&\hfil$#$\hfil&\vrule#&\hfil$#$\hfil&\vrule#&\hfil$#$\hfil&\vrule#&\hfil$#$\hfil&\vrule#&\hfil$#$\hfil&\vrule#&\hfil$#$\hfil&\vrule#&\hfil$#$\hfil&\vrule#&\hfil$#$\hfil&\vrule#&\hfil$#$\hfil&\vrule#&\hfil$#$\hfil&\vrule#\tabskip=0pt\cr
\noalign{\hrule}%
\prored
&k&&1&&2&&3&&4&&5&&6&&7&&8&&9&&10&&11&\cr
\prored
\noalign{\hrule}%
\prored &u_k&&0&&1&&-1&&-2&&9&&-9&&-50&&267&&-413&&-2 180&&17
731&\cr \prored
\noalign{\hrule}%
\prored &v_k&&-1&&1&&1&&-5&&5&&21&&-105&&141&&777&&-5 513&&13
209&\cr \prored
\noalign{\hrule}%
}}
$$

As an illustration of the above summation formulae,  the first
four examples $(\varepsilon = 1)$ are:
$$
\eqalign{&(a) \ \ \ \sum^{n-1}_{i=0} i! \, i = -1+n!\ ,\cr &(b) \
\ \ \sum^{n-1}_{i=0} i! \, (i^2+1) = 1+n! \, (n-1)\ ,\cr &(c) \ \
\ \sum^{n-1}_{i=0} i! \, (i^3-1) = 1+n! \, (n^2-2n-1)\ ,\cr &(d) \
\ \ \sum^{n-1}_{i=0} i! \, (i^4-2) = -5+n! \, (n^3-3n^2+5)\
.\cr}\eqno(14)
$$

Note that $i^k+u_k$ in (11) is a polynomial $P_k(i)$ in (2) in a
reduced form and suitable for generalization. Namely, (11) can be
generalized to
$$
\sum^{n-1}_{i=0} \varepsilon^i \, i! \, P_k(i) =
V_k+\varepsilon^{n-1} n! \, B_{k-1}(n)\ ,\quad k\ge 1\ ,\eqno(15)
$$
where $P_k(i) = \sum_{r=0}^k \, C_r \, i^r$ with
$$ C_0 =
\sum^k_{r=1} \, C_r \, u_r, \quad V_k = \sum^k_{r=1} \, C_r \, v_r
,  \quad B_{k-1}(n) = \sum^k_{r=1} C_r \, A_{r-1}(n)
$$ and
$C_1,\cdots,C_k\in\msbm\hbox{Z}$. Polynomials $P_k(i)$ which do
not have the above form do not yield (15). \vskip1cm

\centerline{{\bf 3. Divisibility}} \vskip.5cm

\noindent The above results enable us to investigate some
divisibility properties of  $\sum^{n-1}_{i=0} \varepsilon^i \, i!
\, P_k(i)$ with respect to all factors contained in $n!$.
According to (15) we have that $\sum^{n-1}_{i=0} \varepsilon^i \,
i! \, P_k(i)$ and $V_k$ are equally divisible with respect to
factors of $n!$, as well as to those of $B_{k-1}(n)$.

\vskip.5cm

\noindent {\bf Proposition 1.} {\it If the polynomial $A_{k-1}(n)$
satisfies the identity (4) then we have the following congruence}
$$
\sum^{n-1}_{i=0} \varepsilon^i \, i!\, [i^k+A_{k-1}(1)-\varepsilon
A_{k-1}(0)] \, \equiv -\varepsilon A_{k-1}(0) \ (\hbox{mod} \
n!).\eqno(16)
$$
\vskip.5cm

{\it Proof}. Congruence (16) is a direct consequence of (9).

From (16) it follows
$$
\sum^{n-1}_{i=0} \varepsilon^i \, i! \, i^k \, \equiv -\,
[A_{k-1}(1)-\varepsilon A_{k-1}(0)]\,
\sum^{n-1}_{i=0}\varepsilon^i \, i! - \varepsilon A_{k-1}(0) \
(\hbox{mod} \ n)
$$
and this property can be used to simplify numerical investigation
of divisibility of $\, \sum^{n-1}_{i=0} \varepsilon^i \, i!\, i^k$
 by $n$.

There is a simple example of (16), {\it e.g.}
$$
\sum^{n-1}_{i=0} i!\, i \, \equiv-1 \ (\hbox{mod} \ n!) \
,\eqno(17)
$$
what follows from  (14.a).

\vskip.5cm

\noindent {\bf Proposition 2.} {\it The following statements are
valid:}
$$
\eqalign{(i) \ \ \ &\sum^{n-1}_{i=0} i! \, i \, \not\equiv0 \
(\hbox{mod} \ n)\ ,\quad n>1\ ,\cr (ii) \ \ \ &\sum^{p-1}_{i=0} i!
\, i \, \not\equiv0 \ (\hbox{mod} \ p)\ ,\quad p\in P\ ,\cr (iii)
\ \ \ &\bigg(\sum^{n-1}_{i=0} i! \, i,n!\bigg) = 1\ ,\quad n>1\
,\cr}\eqno(18)
$$
{\it where $(a,b)$ denotes the greatest common divisor of $a, b\in
\msbm\hbox{Z}$, and $P$ is the set of prime numbers.} \vskip.5cm

{\it Proof.} Every of equations $(i)$, $(ii)$ and $(iii)$ in (18)
follows from  (14.a). One can also show that these statements are
equivalent.

Due to (16) divisibility of $\sum^{n-1}_{i=0} \varepsilon^i \,
i!\, i^k,\ \, k\ge 1$, by factors of $n!$ is in some relation to
divisibility of $\sum^{n-1}_{i=0}  \varepsilon^i \, i!$ except for
the case $A_{k-1}(1) = \varepsilon \, A_{k-1}(0)$.

\vskip1cm \centerline{{\bf 4. On Kurepa's Hypothesis }} \vskip.5cm
Kurepa in [1] introduced a hypothesis
$$
(!n,n!) = 2\ ,\quad 2\le n\in {\msbm \hbox{N}}\ ,\eqno(19)
$$
where
$$
!n = \sum^{n-1}_{i=0}i!\ \eqno(20)
$$
has been called the left factorial. In spite of many papers (for a review
see [2] and references therein) on KH it is still an open problem in
number theory [3].
 Many equivalent statements to KH have been obtained (for some of them see
[4]). Among very simple assertions equivalent to (19) are [1]:
$$
\eqalign{&!n \, \not\equiv0 \ (\hbox{mod} \ n)\ ,\quad n>2\ ,\cr
&!p \, \not\equiv0 \ (\hbox{mod} \ p)\ ,\quad p>2\ .\cr}\eqno(21)
$$
KH is verified by computer calculations (see [2]) for $n<2^{23}$ [5].

The above obtained summation formulae give us possibility to
introduce infinitely many new statements equivalent to KH. The first
three of them, which follow from (14), are:
$$
\eqalign{&\sum^{p-1}_{i=0} i! \, i^2 \, \not\equiv1 \ (\hbox{mod}
\ p)\ ,\quad p>2\ ,\cr &\sum^{p-1}_{i=0} i! \, i^3 \, \not\equiv1
\ (\hbox{mod} \ p)\ ,\quad p>2\ ,\cr &\sum^{p-1}_{i=0} i! \, i^4
\, \not\equiv -5 \ (\hbox{mod} \ p)\ ,\quad p>2\ .\cr}\eqno(22)
$$
\vskip.5cm

\noindent {\bf Theorem 3.} {\it If $u_k$ and $v_k$ satisfy (13.a)
and (13.b) then}
$$
\sum^{p-1}_{i=0} i! \, i^k \, \not\equiv v_k \, (\hbox{mod} \ p)\,
,\quad p>2  \, ,\eqno(23)
$$
{\it is equivalent to KH for such $k\in\msbm\hbox{N}$ for which
$u_k$ is not divisible by $p$.}

\vskip.5cm

{\it Proof.} Consider (11) for $\varepsilon=1$ and $n=p$.
According to KH one has $ u_k\, \sum^{p-1}_{i=0} i! \,\not\equiv 0
\, (\hbox{mod} \ p) $ for $p>2$ and $p$ which does not divide
$u_k\not=0$. For such primes $p$ it holds (23).

Starting from the  Fermat little theorem, i.e. $i^{p-1}=1$ in the
Galois field GF($p$) if $i=1,2,...,p-1$\ , one can easily show
that assertion
$$
\sum^{p-1}_{i=0} i! \, i^{r(p-1)} \, \not\equiv -1 \ (\hbox{mod} \
p)\ ,\quad p>2\ , \ r\in\msbm\hbox{N}           \eqno(24)
$$
is equivalent to KH. This can be regarded as a special case of the
Theorem 3. Since $r$ may be any positive integer it means that
there are infinitely many equivalents to KH.

 Note that on the basis of  Fermat's theorem  one can also obtain
\vskip.5cm
$$
\sum^{p-1}_{i=0} \varepsilon^i \, i! \, i^{k+r(p-1)} =
\sum^{p-1}_{i=0} \varepsilon^i \, i!\, i^k - \delta_{0k}\ ,\quad
k\in {\msbm\hbox{N}}_0\ ,\quad   r\in \msbm\hbox{N}\ .\eqno(25)
$$

\vskip.5cm

\noindent {\bf Proposition 3.}  { \it If $u_k$ and $v_k$ satisfy
(13.a) and (13.b), respectively,  the following relations in
GF$(p)$ are valid:}
$$
(u_{p-1}+1)\, \sum^{p-1}_{i=0} i! = v_{p-1}+1\ ,\eqno(26.a)
$$
$$
u_p \, \sum^{p-1}_{i=0} i! = v_p+1\ ,\eqno(26.b)
$$
$$
(u_{p+1}-1)\sum^{p-1}_{i=0}i! = v_{p+1}-1\ ,\eqno(26.c)
$$
$$
(u_{p+2}+1) \, \sum^{p-1}_{i=0}i! = v_{p+2}-1\ .\eqno(26.d)
$$
\vskip.5cm

{\it Proof.}  One can  start from (11), then use (25) and (14).

From eqs. (13.a) and (13.b)  one obtains in GF($p$):

$(u_{p+2},v_{p+2})=(-u_p-1,-v_p)$,

$ (u_{p+1},v_{p+1})=(1,1)$,

$(u_p,v_p)=(u_{p-1}+1,v_{p-1}) .$

\noindent Thus $(26.a)-(26.d)$  are equivalent identities which
are always satisfied owing to the values of $u_k$ and $v_k$.
Identity (26.c) does  not depend on validity of KH.

\vskip1cm
\centerline{{\bf 5. Concluding remarks}}
\vskip.5cm

It is worth noting that for every $k\in \msbm\hbox{N}$ there is a
unique pair $ (u_k,v_k) $ of integers $u_k$ and $v_k$ which
connect $ \sum^{n-1}_{i=0} \varepsilon^i \, i!\, i^k $ and $
\sum^{n-1}_{i=0} \varepsilon^i \, i! $ into simple summation
formula (11). All other results of the present paper are mainly
various consequences of this fact.

Formula (11) is also suitable to consider its limit when $ n\to \infty
$ in $p$-adic analysis. Namely, since $ |n!|_p \to 0 $ as
$ n\to \infty $, one obtains
$$
\sum^\infty_{i=0} \varepsilon^i \, i! \,(i^k+u_k)=v_k \ ,
$$
valid in $ \msbm \hbox{Q}_p $ for every $p$. Some p-adic aspects
of the series $\sum^\infty_{i=0} \varepsilon^i \, i! \, P_k(i)$
and their possible role in theoretical physics are considered in
Ref. 6.

Having infinitely many new equivalents, Kurepa's hypothesis
becomes more challenging. Moreover, KH itself seems to be the
simplest among all its equivalents. In $p$-adic case KH can be
also formulated as follows:
$$
\sum^\infty_{i=0}i! = a_0 + a_1p + a_2p^2 +... , \quad  p\in P \ ,
$$
where $a_i$ are definite digits with $ a_0\ne 0 $ for all $ p\ne 2
$.

\vskip1cm
\centerline{{\bf Acknowledgment}}
\vskip.5cm

The author thanks  \v Z. Mijajlovi\'c for useful discussions.
\vskip1cm
\centerline{References}
\vskip.5cm

\item{[1]} Dj. Kurepa, {\it On the left factorial function}, Math.
Balkanica  {\bf 1} (1971), 147-153. \item{[2]} A. Ivi\'c and \v Z.
Mijajlovi\'c, {\it On Kurepa's problems in number theory}, Publ.
Inst. Math {\bf 57} (71) (1995), 19-28. \item{[3]} R. Guy, {\it
Unsolved Problems in Number Theory}, Springer-Verlag, 1981.
\item{[4]} Z.N. \v Sami, {\it A sequence $u_{n,m}$ and Kurepa's
hypothesis on left factorial}, Scientific Review {\bf 19-20}
(1996), 105-113. \item{[5]} M. \v Zivkovi\'c, {\it On Kurepa left
factorial hypothesis}, Kurepa's Symposium, Belgrade 1996.
\item{[6]} B. Dragovich, {\it On  some p-adic series with
factorials}, Lecture Notes in Pure and Applied Mathematics {\bf
192} (1997), 95-105 ; math-ph/0402050.

\vskip1.5cm

Institute of Physics

P.O.Box 57 , 11001 Belgrade

Yugoslavia

dragovich@phy.bg.ac.yu
\end